\documentclass[a4paper,11pt,reqno]{amsart}

\newcommand{\erre}{\mathbb{R}}
\renewcommand{\r}{\mathbb{R}}
\newcommand{\n}{\mathbb{N}}
\newcommand{\E}{\mathbb{E}}
\newcommand{\EE}{\mathcal{E}}
\renewcommand{\P}{\mathbb{P}}

\newcommand{\ds}{\displaystyle}
\newcommand{\osc}{\mathop{\mathrm{osc}}\nolimits}

\newcommand{\var}{\mathop{\mathrm{Var}}\nolimits}

\newtheorem{prop}{Proposition}
\newtheorem{thm}[prop]{Theorem}
\newtheorem{coroll}[prop]{Corollary}

\theoremstyle{definition}
\newtheorem{rmk}[prop]{Remark}

\begin{document}

\title[Convergence of MCMC methods]{Convergence of sequential Markov Chain Monte Carlo methods: I. Nonlinear flow of probability measures}
\author{Andreas Eberle}
\address{Institut f\"ur Angewandte Mathematik, Universit\"at Bonn, Wegelerstr. 6, 53115 Bonn, Germany}
\email{eberle@uni-bonn.de}
\urladdr{http://wiener.iam.uni-bonn.de/$\sim$eberle}

\author{Carlo Marinelli}
\address{Institut f\"ur Angewandte Mathematik, Universit\"at Bonn, Wegelerstr. 6, 53115 Bonn, Germany}
\urladdr{http://wiener.iam.uni-bonn.de/$\sim$marinelli}

\date{December 3, 2006}

\begin{abstract}
Sequential Monte Carlo Samplers are a class of stochastic
algorithms for Monte Carlo integral estimation w.r.t.\ probability
distributions, which combine elements of Markov chain Monte Carlo
methods and importance sampling/resampling schemes. We develop a
stability analysis by functional inequalities for a nonlinear flow
of probability measures describing the limit behavior of the
algorithms as the number of particles tends to infinity. Stability
results are derived both under global and local assumptions on the
generator of the underlying Metropolis dynamics. This allows us to
prove that the combined methods sometimes have good asymptotic
stability properties in multimodal setups where traditional MCMC
methods mix extremely slowly. For example, this holds for the mean
field Ising model at all temperatures.
\end{abstract}

\subjclass[2000]{65C05, 60J25, 60B10, 47H20, 47D08}

\keywords{Markov Chain Monte Carlo, sequential Monte Carlo,
importance sampling, spectral gap, Dirichlet forms, functional
inequalities.}

\thanks{This work is supported by the Sonderforschungsbereich 611,
  Bonn. The second author also gratefully acknowledges the hospitality
  and financial support of the Institute of Mathematics, Polish
  Academy of Sciences, Warsaw, of the Institute des Hautes \'Etudes
  Scientifiques, Bures-sur-Yvette, of the Max-Planck-Institut f\"ur
  Mathematik, Leipzig, through an IPDE fellowship.}

\maketitle

Spectral gap estimates, or, equivalently, Poincar\'e inequalities,
as well as other related functional inequalities provide powerful
tools for the study of convergence to equilibrium of reversible
time-homogeneous Markov processes (see e.g. \cite{DSC-comp},
\cite{DSC-logsob}, \cite{DSC-nash}). In particular, they have been
successfully applied to analyze convergence properties of Markov
Chain Monte Carlo (MCMC) methods based on reversible Markov chains
(see e.g.\ \cite{DSC-metro}). The idea of MCMC methods is to
produce approximate samples from a probability distribution $\mu$
by simulating for a sufficiently long time an ergodic Markov chain
having $\mu$ as invariant measure. MCMC methods have become the
standard to carry out Monte Carlo integrations with respect to
complex probability distributions in many fields of applications,
including in particular Bayesian statistics, statistical physics,
and computational chemistry. We just refer the interested reader to \cite{liu}
and \cite{casella} and references therein as an example of work in
this area, as the literature is by now enormous.
Since the Markov chain is usually started with an initial
distribution that is very different from $\mu$, strong convergence
properties, such as exponential convergence to equilibrium with a
sufficiently large rate, are required to ensure that the
corresponding MCMC method produces sufficiently good approximate
samples from $\mu$. However, these strong convergence properties
often do not hold in multimodal, and in particular
high-dimensional problems, as they arise in many applications. For
example, in statistical mechanics models with phase transitions,
the rate of convergence often decays exponentially in the system
size within the multi-phase regime.\smallskip

In this and a follow-up article we initiate a study of convergence
properties by functional inequalities for a class of algorithms
for Monte Carlo integral estimation that are a combination of
sequential Monte Carlo and MCMC methods. Instead of producing
constantly improved samples of a fixed distribution $\mu$, these
{\em sequential MCMC methods} try to keep track as precisely as
possible of an evolving sequence $(\mu_t)_{0\leq t \leq
  \beta }$ of probability distributions. Here $\mu_0$ is an initial distribution
that is easy to simulate, and $\mu_\beta $ is the target
distribution that one would like to simulate. Importance sampling
and resampling steps are included to constantly adjust for the
change of the underlying measure. Whereas for MCMC methods
exponential asymptotic stability is usually required to obtain
improved samples, the sequential MCMC method starts with a good
estimate of $\mu_0$, and one only has to control the growth of the
``size'' of the error. As a consequence, the method sometimes
works surprisingly well in multimodal situations where traditional
MCMC methods fail, cf.\ also the examples below. The price one has
to pay is that samples from $\mu_\beta $ cannot be produced
individually. Instead, the corresponding algorithm produces
directly a Monte Carlo estimator for $\mu_\beta $ given by the
empirical distribution of a system of interacting particles at the
final time. To ensure good approximation properties, a large
number $N$ of particles is required.\smallskip

Variants of such sequential MCMC methods have recently been proposed
at several places in the statistics literature, see in particular
\cite{DMDJ} and references therein, as well as \cite{cappe}.  However,
precise and general mathematical methods for the convergence and
stability analysis, in the spirit of those developed for traditional
MCMC methods by Diaconis, Saloff-Coste, Jerrum, Sinclair, and many
others, seem still to be missing -- although very important first
steps can be found in the work of Del Moral and coauthors, cf.\ e.g.\
\cite{DMD}, \cite{DMG} and \cite{MM}.  The classical approach via
Dobrushin contraction coefficients is usually limited to very regular
situations.  Moreover, it rarely yields precise statements on the
convergence properties, and it can not be combined easily with
decomposition techniques.\smallskip

Our aim is to make variants of the powerful techniques of the
spectral gap/Dirichlet form approach to convergence rates of
time--homogeneous Markov chains (e.g.\ canonical paths, comparison
and decomposition results) available in the different context of
sequential MCMC methods. Mathematically, this means at first to
study a class of nonlinear evolutions of probability measures by
functional inequalities. Such a study has been initiated in a
related context by Stannat \cite{St-var}. In this work, we
restrict ourselves to the simplest and most natural variant of
sequential MCMC, where importance sampling/resampling is
\emph{only} used to adjust constantly for the change of the
underlying distribution, and MCMC steps at time $t$ are always
carried out such that detailed balance holds w.r.t.\ the measure
$\mu_t$ (and not w.r.t.\ $\mu_0$!). This seems crucial for
establishing good stability properties. Note that the type of
sequential Monte Carlo samplers studied here is different from
those analyzed by Del Moral and Doucet in \cite{DMD}. An
algorithmic realization has been applied to simulations in
Bayesian mixture models by Del Moral, Doucet and Jasra in
\cite{DMDJ}, who observed substantial benefits compared to other
methods.
\smallskip

We have divided our work on sequential MCMC methods into two
publications: in this first article we study the stability properties
of nonlinear flows of probability measures describing the limit as the
number $N$ of particles goes to infinity. In the follow-up work
\cite{EM2} we will apply the results to control the asymptotic
variances of the Monte Carlo estimators as $N$ tends to infinity. The
functional inequality approach enables us to prove stability
properties not only under global but also under local conditions, i.e.
assuming only that good estimates hold on each set of a decomposition
of the state space. As a consequence, we obtain a procedure for
analyzing the asymptotic behavior of sequential MCMC methods applied
to multimodal distributions. For example, in the spirit of previous
results for tempering algorithms by Madras and Zheng \cite{MZ} and
others, we can prove good (polynomial in the system size and the
inverse temperature) stability properties in the case of the mean
field Ising model, cf. Section \ref{subsection:ising} below. We also
demonstrate in a simple exponential model with several modes that the
flow of probability measures corresponding to sequential MCMC methods
has better stability properties than the one corresponding to the
classical simulated annealing algorithm, cf.
\ref{subsection:expmodel}. Sequential MCMC methods might hence also
provide an efficient alternative to simulated annealing.

\section{Setup}

\subsection{Sequential estimation of probability measures}

Let $S$ denote a finite state space, and $\mu$ a probability
measure on $S$ with full support, i.e. $\mu (x)>0$ for all $x\in
S$. The finiteness of the state space is only assumed to keep the
presentation as simple and non-technical as possible. Most results
of this paper extend to continuous state spaces under standard
regularity assumptions. By ${\mathcal M}_1(S)$ we denote the space
of probability measures on $S$. As usual,
$$
\nu (f)\ :=\ \int_S f\,d\nu\ =\ \sum_{x\in S} f(x)\, \nu (x)
$$
denotes the expectation of a function $f:S\to\r $ w.r.t.\ a
measure $\nu\in {\mathcal M}_1(S)$. We consider methods for Monte
Carlo integration with respect to the probability distributions of
an exponential family
\begin{equation}
  \label{eq:mub}
  \mu_t(x) = \frac1{Z_t}\, e^{-t H(x)}\mu(x),
             \quad 0 \leq t < \infty,
\end{equation}
where $H:S\to [0,\infty )$ is a given function, and
$Z_t:=\sum_{x\in S} e^{-t H(x)}\,\mu(x)$ is a normalization constant.
Below, $t$ will play the r\^{o}le of a time parameter for a particle
system approximation.

Note that for a fixed $\beta >0$, {\em any} given probability
measure $\nu$ on $S$ that is mutually absolutely continuous with
respect to $\mu$ can be written in the form (\ref{eq:mub}) with
$t=\beta $ by setting $H(x)=\frac1\beta
\log\frac{\mu(x)}{\nu(x)}$. One should then think of the family
$(\mu_t)_{0\leq t\leq \beta }$ of probability measures as a
particular way to interpolate between the target distribution
$\mu_{\beta }$ that we would like to simulate, and the reference
distribution $\mu_0=\mu $ that can be simulated more easily.
Although we restrict our attention here to this simple way of
interpolating between two measures, other interpolations can be
treated by similar methods. In fact, an arbitrary family
$\left(\mu_t\right)_{0\leq t\leq\beta }$ of mutually absolutely
continuous probability measures on $S$ with smooth dependence on
$t$ can be written in the form
 \begin{equation}
 \label{eq:timedep}
 \mu_t(x)\ =\ \frac 1{Z_t}\,e^{-\int_0^tU_s(x)\,ds}\,\mu (x)\,
 ,\qquad 0\leq t\leq\beta\, ,
 \end{equation}
 where $Z_t$ is a normalization constant, and $(s,x)\mapsto
 U_s(x)$ is a continuous non--negative function on $[0,\beta ]\times S$.
 Our results below extend to this more general case, cf.\ Section
 \ref{subsection:extensions} below.

  The main advantage of the interpolation
(\ref{eq:mub}) is that the singularity of $\mu_t$ w.r.t.\ $\mu$ is
resolved uniformly over time. In particular,
\begin{equation}
  \label{eq:rd}
  \left|\log\frac{\mu_{s}(x)}{\mu_t (x)}\right| \ \leq\
  \osc(H)\cdot |s-t|\qquad \mbox{for all } s,\, t\geq 0\mbox{ and
  }x\in S
\end{equation}
where $\osc(H) := \max_S H - \min_S H$. On the other hand, other
interpolations, e.g.\ by a spatial coarse graining, may be
preferable in concrete applications.\medskip

One way to obtain sequential methods for Monte Carlo estimation of
expectation values with respect to the measures $\mu_t$ is to
proceed as follows:
\begin{enumerate}
\item[a)] Construct a semigroup $(\Phi_{s,t})_{0\leq s\leq
    t<\infty}$ of nonlinear transformations on the space of
  probability measures on $S$, such that
\begin{equation}
  \label{eq:flow}
\Phi_{s,t}\mu_s = \mu_t \qquad \mbox{ for all } 0 \leq s \leq t\,
.
\end{equation}
\item[b)] {\em Spatial discretization by interacting particle
system}: Construct an appropriate Markov process $(X_t^1,\ldots
,X_t^N)$ on $S^N$ ($N\in\n $) related to the nonlinear semigroup
$\Phi_{s,t}$, and estimate $\mu_t=\Phi_{0,t}\mu$ by the empirical
distributions
$$
\hat\mu_t^{(N)}\ :=\ \frac1N\, \sum_{i=1}^N \delta_{X_t^i}, \qquad
t\geq 0,
$$
of the process with initial distribution $\mu^N$.
\item[c)] {\em Time--discretization}: Approximate the continuous time
  Markov process $(X_t^1,\ldots,X_t^N)$ by a time--discrete Markov
  chain on $S^N$ (which can then be simulated).
\end{enumerate}

\subsection{The non-linear semigroup}

To define the nonlinear semigroup $\Phi_{s,t}$ and the particle
system we have in mind, we consider the generators ($Q$-matrices)
${\mathcal L}_t$ at time $t\ge 0$ of a time-inhomogeneous Markov
chain on $S$ satisfying the detailed balance condition
\begin{equation}
  \label{eq:db}
  \mu_t(x){\mathcal L}_t(x,y) = \mu_t(y){\mathcal L}_t(y,x)
      \quad \forall\ t\geq0, \; x,y\in S.
\end{equation}
The generators ${\mathcal L}_t$ determine the MCMC steps in a
corresponding sequential MCMC method. We assume that ${\mathcal
L}_t(x,y)$ depends continuously on $t$. To compare algorithmic
performance in a reasonable way, one might also assume
\begin{equation}
 \label{eq:star}
\sum_{y\neq x}{\mathcal L}_t(x,y)\leq 1\qquad\forall\ x\in S\, ,
\end{equation}
although this is not necessary for the results below. For example,
${\mathcal L}_t$ could be the generator of a Metropolis dynamics
w.r.t.\ $\mu_t$, i.e.,
$${\mathcal L}_t(x,y)\ =\ K_t(x,y)\cdot\min\left(
\frac{\mu_t(y)}{\mu_t (x)} ,1\right)\mbox{ for }x\neq y,$$
${\mathcal L }_t(x,x)=-\sum_{y\neq x}{\mathcal L}_t(x,y)$, where
the proposal matrix $K_t$ is a given symmetric transition matrix
on $S$. By (\ref{eq:db}), ${\mathcal L}_t$ defines a symmetric
linear operator on $L^2(S,\mu_t)$. The associated Dirichlet form
on functions $f,g:S\to\r$ is
$$
\EE_t(f,g)\ :=\ -\E_t[f\,{\mathcal L}_tg]\ =\ \frac
12\,\sum_{x,y\in S}(f(y)-f(x))(g(y)-g(x))\,{\mathcal
L}_t(x,y)\,\mu_t(x)\, ,
$$ where $\E_t$ stands for expectation w.r.t.\ $\mu_t$, and
$$({\mathcal L}_tg)(x)\ :=\ \sum_y{\mathcal L}_t(x,y)g(y)\, .$$
We shall often
use the abbreviated notation $\EE_t(f):=\EE_t(f,f)$.

We fix non--negative constants $M_t$ ($t\ge 0$) that determine the
average relative frequency of MCMC moves compared to importance
sampling/resampling steps in a corresponding SMCMC method. Again,
we assume that $t\mapsto M_t$ is continuous.\smallskip

Let $p_{s,t}(x,y)$ and $q_{s,t}(x,y)$ ($x,y\in S$) be the unique
solutions of the forward equations
\begin{eqnarray}
\label{eq:fp} \ds \frac{\partial}{\partial t}\, p_{s,t} f &=&
p_{s,t}(M_t{\mathcal L}_tf-Hf),
      \qquad p_{s,s}f=f,\\
 \label{eq:fq}
\ds \frac{\partial}{\partial t}\, q_{s,t} f &=&
q_{s,t}(M_t{\mathcal L}_tf-H_tf),
      \qquad q_{s,s}f=f,
\end{eqnarray}
where $$H_t\ :=\ H-\E_t[H]\ .$$
The linear semigroups $p_{s,t}$ and $q_{s,t}$ have the Feynman-Kac
representations
\begin{eqnarray*}
  p_{s,t}f(x) &=& \int
                  e^{-\int_s^t H(X_r(\omega ))\,dr}f(X_t(\omega ))\; \P_{t,x}(d\omega )\, ,\\
  q_{s,t}f(x) &=& \int
                  e^{-\int_s^t H_r(X_r(\omega ))\,dr}f(X_t(\omega ))\; \P_{t,x}(d\omega )\, ,
\end{eqnarray*}
where $(X_t,\P_{t,x})$ is a time-inhomogeneous Markov process with
generator $M_t\cdot {\mathcal L}_t$. In particular one has
$$
q_{s,t}f = \exp\Big( \int_s^t \E_r[H]\,dr \Big)\, p_{s,t}f.
$$
We consider the nonlinear semigroup
$$
\Phi_{s,t}\nu\ := \nu\ \frac{\nu p_{s,t}}{(\nu p_{s,t})(S)}\ =\
\frac{\nu q_{s,t}}{(\nu q_{s,t})(S)}, \qquad 0 \leq s \leq t,
$$
on the space ${\mathcal M}_1(S)$ of probability measures on $S$.
Here
$$
\nu p(y) = \sum_{x\in S}\nu(x)p(x,y) \, .
$$
The semigroup $\Phi_{s,t}$ describes the time
evolution of the law of an inhomogeneous Markov chain with
generator $M_t\cdot{\mathcal L}_t$ and absorption rate $H$,
conditioned to be alive at time $t$ (see e.g. \cite{BG}). It is
not difficult to verify that (\ref{eq:flow}) holds, cf.\ Theorem
\ref{thm:err} below.

\subsection{Particle system approximations}

To approximate $\Phi_{s,t}$ one could use a particle system
consisting of independent Markov chains with absorption, and base
the Monte Carlo estimation on the particles that are still alive
at time $t$. However, such a procedure would be usually very
inefficient, since in most interesting cases the overwhelming
majority of particles would have become extinct already at the
final time. Instead, sequential Monte Carlo samplers are based on
a time--inhomogeneous Markov chain on $S^N$, $N\in\n $, with a
generator that is for example of type
\begin{multline*}
\bar{\mathcal L}_t^Nf\, (x^1,\ldots ,x^N)\ =\ M_t\cdot
\sum_{i=1}^N({\mathcal L}_t^{(i)}f)(x^1,\ldots ,x^N) \\
 + \frac1N \sum_{i,j=1}^NH(x^i)\cdot (f(r^{i,j}(x))-f(x))\, .
\end{multline*}
Here ${\mathcal L}_t^{(i)}$ denotes the application of ${\mathcal
L}_t$ to the $i$-th component, and $r^{i,j}(x):=y$ where
$y^i:=x^j$ and $y^k:=x^k$ for all $k\neq i$. Hence, between the
interactions the particles move according to time--inhomogeneous
Markov chains with generator $M_t\cdot {\mathcal L}_t$ and
absorption rate $H$, and in case of absorption, the position is
replaced by the position of a randomly chosen particle. Other
interaction terms that correspond to different resampling schemes
are possible as well. The asymptotics as $N\to\infty$ of the
approximating particle systems with mean field interaction has
been studied intensively, cf.\ e.g.\ the monograph \cite{DM}.

\subsection{Convergence and stability properties}

The quality of Monte Carlo estimates of $\mu_t(f)=\int f\, d\mu_t$
for some function $f:S\to\erre$ can be measured by the bias and
the (asymptotic) variance of the corresponding estimators. The
theoretical analysis of the sequential MCMC methods considered
here can be subdivided into several steps as above:
\begin{enumerate}
\item[a)] Stability properties of the semigroup $\Phi_{s,t}$.
\item[b)] Bias and asymptotic variance of the estimators
$\hat\mu_t^N(f)=\frac 1N\sum_{i=1}^Nf(X_t^i)$.
\item[c)] Effect of the discretization in time.
\end{enumerate}

In this paper, we will focus exclusively on the first step, that
is we develop a stability analysis for $\Phi_{s,t}$ based on
functional inequalities. A follow-up paper \cite{EM2} will be
devoted to the time dependence of the asymptotic (as $N\to\infty$)
mean square error of the particle system based estimators
$\hat\mu_t^N(f)$. Let us remark for the moment, that significant
work in this direction has already been done, e.g., by Del Moral
and Miclo in \cite{MM}. The results clearly indicate that
techniques very close to those developed here can also be applied
to control the asymptotic variances of the approximating particle
systems. This will be made precise in \cite{EM2}.

We also point out that usually the time discretization is carried
out before the spatial discretization, i.e. one usually directly
considers semigroups and particle systems in discrete time. Even
though this is closer to the algorithmic realization, the
convergence analysis becomes more transparent in continuous time
due to the infinitesimal description (at least from an analytic
perspective). Moreover, the continuous time setup allows us to see
more clearly how frequently different types of moves of the
particle systems should be carried out.

Before stating our results, we comment on relations of sequential
MCMC methods to several standard methods for Monte Carlo
integration:

--~\emph{Parallel MCMC} is a special case of the algorithm above
when $H\equiv 0$, i.e. $\mu_t = \mu$ for all $t$. In this case the
associated particle system consists of independent {\em
time-homogeneous} Markov chains with invariant measure $\mu$.
Common problems are slow mixing due to multimodality and the
burn-in time (i.e. the time needed to reach equilibrium from an
initial distribution that is far from $\mu$ can be much larger
than the inverse spectral gap). Both problems are particularly
significant in high dimensional setups (``curse of dimension'').

--~In \emph{parallel simulated annealing}, the approximating
particle system is given by independent time--inhomogeneous Markov
chains with generator ${\mathcal L}_t$. There are no
interactions. In this case, the corresponding (linear) semigroup on
probability measures does not satisfy (\ref{eq:flow}). As a
consequence, there is an asymptotic bias of the corresponding
Monte Carlo estimator, which can only be reduced by the mixing
properties of the underlying Markov chains. Therefore, in
multimodal setups good convergence properties can only be
guaranteed if the measures $\mu_t$ change very slowly (logarithmic
cooling schedule).

--~Pure \emph{importance sampling/resampling} is the special case
of our method when ${\mathcal L}_t\equiv 0$ for all $t$. Since the
particles cannot explore the state space, it is only applicable
for small state spaces, or in very special situations. In fact,
the results below indicate that a certain amount of particle
motion is needed to ensure good stability properties. Our results
below can be used to quantify, at least in principle, how many MCMC
moves are needed to balance the error growth due to importance
sampling/resampling.

--~A combination of importance sampling and MCMC (without
resampling) is similar to considering Markov chains with
absorption, conditioned to stay alive. This is often inefficient,
cf. the remark above.

--~Finally, we would like to point out that the analysis of several
multilevel sampling methods (see e.g. \cite{liu}) such as
\emph{umbrella sampling} (cf. \cite{MP}), \emph{simulated} and
\emph{parallel tempering} (cf. \cite{MZ}, \cite{BR}, \cite{sch}) has
been an inspiration for this work. These MCMC methods provide samples
from mixtures, direct sums, or products of distributions $\mu_{t_i}$
($0\le i\le m$), $0=t_0<t_1<\cdots <t_m$, $m\in\n $. A disadvantage of
umbrella sampling and simulated tempering is that the normalization
constants have to be estimated in parallel.  Parallel tempering avoids
this disadvantage by simulating the product distribution
$\bigotimes_{i=0}^m\mu_{t_i}$ with the help of a Metropolis chain on
$S^{m+1}$ where neighboring coordinates can be swapped.  However, the
swapping procedure seems to slow down the convergence in some cases,
and it makes the convergence analysis rather intricate, cf. \cite{MZ}.
The sequential MCMC methods presented here can be seen as an attempt
to overcome these difficulties. The estimation of the normalization
constant is built into the algorithm, and the evolution in $t$ is
linear -- and thus faster than a diffusive motion in $t$ as in
simulated and parallel tempering. Once the basic techniques are
developed, the convergence analysis seems also to be at least
partially more transparent for sequential MCMC than for simulated and
parallel tempering.

\section{Main results}

\subsection{Time evolution of the mean square error}

Let $\nu_t:=\Phi_{0,t}\nu$ for some given initial distribution
$\nu\in {\mathcal M}_1(S)$, and let
$$
g_t(y) := \frac{\nu_t(y)}{\mu_t(y)}, \qquad t\ge 0,
$$
denote the relative density of $\nu_{t}$ w.r.t.\ the measure
$\mu_t$ defined by (\ref{eq:mub}). Moreover, let
$$
{\varepsilon}_t := \E_t[(g_t-1)^2]
$$
denote the mean square error ($\chi^2$--contrast) of $\nu_t$
w.r.t.\ $\mu_t$. Our first result shows that
$\mu_t=\Phi_{0,t}\mu_0$, and it gives a general method to analyze
the stability of this evolution in an $L^2$ sense:

\begin{thm}\label{thm:err}
(i) $\nu_t=\Phi_{0,t}\nu$ is the unique solution of the nonlinear
evolution equation
\begin{equation}
\label{eq:nut} \frac{\partial}{\partial t} {\nu}_t\  =\
M_t\nu_t{\mathcal L}_t\, -\, H\,\nu_t\, +\, \nu_t(H)\, \nu_t,
\qquad t\ge 0.
\end{equation}
with initial condition $\nu_0=\nu$.
\medskip

(ii) The densities $g_t$ solve
\begin{equation}
\label{eq:X}
\frac{\partial}{\partial t} {g}_t\  =\ M_t{\mathcal
L}_tg_t\, +\, \E_t[H (g_t-1)]\, g_t.
\end{equation}
(iii) The time evolution of the mean square error is given by
\begin{equation}
\label{eq:err} \frac12\frac{d}{dt} {\varepsilon}_t =
-M_t\,\mathcal{E}_t(g_t-1) - \frac12 \E_t[H_t (g_t-1)^2] +
\E_t[H_t(g_t-1)]\, {\varepsilon}_t.
\end{equation}
\end{thm}\smallskip

\begin{rmk}
(\ref{eq:nut}) is the forward equation for the nonlinear semigroup
$\Phi_{s,t}$ (for $s=0$). The corresponding assertion holds for
$\nu_t:=\Phi_{s,t}\nu $ for $t\geq s>0$. Since $\mu_t$ solves
(\ref{eq:nut}), we obtain in particular
$$\mu_t=\Phi_{s,t}\mu_s\qquad\mbox{for all }t\geq s\geq 0\, .$$
\end{rmk}

The proof of the theorem is given in Section \ref{sec:proofs} below.
Similar equations have been derived in a more general setup by Stannat
\cite{St-var}.
\medskip

The main objective of this article is to develop efficient tools
to bound the growth of $\varepsilon_t$ based on Theorem
\ref{thm:err}. To estimate the right-hand side of (\ref{eq:err})
we have to control the two terms involving $H_t$ (which correspond
to importance sampling/resampling) by the Dirichlet form
$\mathcal{E}_t$ (which corresponds to MCMC moves). We first
discuss how this can be achieved in the presence of a good global
spectral gap estimate. Afterwards, we give results based on local
Poincar\'{e}-type inequalities, which can sometimes be used to
control the error growth in multimodal setups where good global
mixing properties of the underlying Markov chains do not hold.

\subsection{Stability based on global estimates}
\label{subsection:global} For $t\ge 0$ let
$$C_t\ :=\ \sup\
\left\{\left.\E_t[f^2]/\mathcal{E}_t(f,f)\,\right|\;
f:S\to\r\mbox{ s.t.\ }\E_t[f]=0\, ,\; f\not\equiv 0\,\right\}
$$
denote the (possibly infinite) inverse spectral gap of
$\mathcal{L}_t$, and let
$$A_t\ :=\ \sup\
\left\{\left.\E_t[H_t^-f^2]/\mathcal{E}_t(f,f)\,\right|\;
f:S\to\r\mbox{ s.t.\ }\E_t[f]=0\, ,\; f\not\equiv 0\,\right\}\, .
$$
Thus $C_t$ and $A_t$ are the optimal constants in the global
Poincar\'{e} inequalities
\begin{eqnarray}
\label{eq:Poin} {\rm Var}_t(f)& \leq & C_t\cdot \mathcal{E}_t(f,f)
\qquad\forall\ f:S\to\r\, ,
\qquad \mbox{and}\\
  \label{eq:H}
\E_t[H_t^- (f-\E_t[f])^2] &\leq & A_t\cdot
\mathcal{E}_t(f,f)\qquad\forall\ f:S\to\r \, .
\end{eqnarray}
Here ${\rm Var}_t$ denotes the variance w.r.t.\ $\mu_t$.

\begin{rmk}
(i) There exist efficient techniques to obtain upper bounds for
$C_t$, for example the method of canonical paths, comparison
methods (see e.g. \cite{SC}), as well as decomposition methods
(see e.g. \cite{JSTV}). Variants of these techniques can be
applied to estimate $A_t$ as well.

(ii) Clearly, one has
\begin{equation}
\label{eq:atct}
A_t\ \le \ C_t
\cdot\sup_{x\in S}{H}_t^-(x),
\end{equation}
so an upper bound on $C_t$ yields a trivial (and usually far from
optimal) upper bound on $A_t$.
\end{rmk}

Let
$$
\sigma_t(H)\ :=\ \var_t (H)^{1/2}\ =\ \E_t[H_t^2]^{1/2}
$$
denote the standard deviation of $H$ w.r.t.\ $\mu_t$. The next
result bounds the error growth in terms of $C_t$ and $A_t$~:

\begin{thm}
\label{thm:logerr} If $M_t\geq  A_t/2$ for all $t\geq 0$, then
\begin{eqnarray}
\label{eq:log1} \frac{d}{dt}\, \log {\varepsilon}_t \leq
-\frac{2M_t-A_t}{C_t}\, +\, 2\,\sigma_t(H) {\varepsilon}_t^{1/2}
\end{eqnarray}
and
\begin{eqnarray}
\label{eq:log2} \frac{d}{dt}\, \log {\varepsilon}_t \leq
-\frac{2M_t-A_t}{C_t} + 2
\Big(\frac{A_t}{C_t}\E_t[{H}_t^-]\Big)^{1/2} {\varepsilon}_t^{1/2}
+ \E_t[{H}_t^-]\,{\varepsilon}_t.
\end{eqnarray}
\end{thm}\smallskip

The proof will be given in Section \ref{sec:proofs} below. Inequality
(\ref{eq:log1}) is straightforward to prove, but sometimes
(\ref{eq:log2}) is stronger, since the constants only depend on
the negative part of $H_t$.
As an immediate consequence of the theorem we obtain estimates on
the average relative frequency $M_t$ of MCMC moves that is
sufficient to guarantee stability of the corresponding nonlinear
flow of probability measures:

\begin{coroll}\label{coroll:glob}
Let $0\leq\beta_0<\beta_1$, and assume that for all $t\in
(\beta_0,\beta_1)$,
\begin{equation}
\label{eq:steps1} M_t \ >\ \frac{A_t}2\, +\, C_t
\sigma_t(H)\,\varepsilon_{\beta_0}^{1/2}
\end{equation}
or
\begin{equation}
\label{eq:steps2} M_t \ >\  \frac{A_t}2\, +\, (A_tC_t
\E_t[H_t^-])^{1/2}\, \varepsilon_{\beta_0}^{1/2} + \frac 12\, C_t
\E_t[H_t^-] \, \varepsilon_{\beta_0}\ .
\end{equation}
Then $t\mapsto {\varepsilon}_t$ is strictly decreasing on the
interval $[\beta_0,\beta_1]$.
\end{coroll}\smallskip

\begin{rmk}
(i) On the finite state spaces considered here, the constants
$C_t$ and $A_t$ are finite if ${\mathcal L}_t$ is irreducible.
However, in multimodal situations, the numerical values of these
constants are often extremely large. Alternative estimates based
on local Poincar\'{e}-type inequalities are given below.

(ii) Similarly to the corollary, one obtains that the error decays
exponentially with rate $\gamma >0$, i.e. $t\mapsto e^{\gamma
t}\,{\varepsilon}_t$ is decreasing on $[\beta_0,\beta_1]$,
provided
\begin{equation}
\label{eq:expdecay} M_t\ >\ \frac{A_t+\gamma C_t}2 \, + C_t
\sigma_t(H)\, e^{-\gamma(t-\beta_0)/2} \,
\varepsilon_{\beta_0}^{1/2}\qquad\forall\ t\in (\beta_0,\beta_1)\,
,
\end{equation}
or a similar
condition replacing (\ref{eq:steps2}) holds.

(iii) One can often assume that the initial error
$\varepsilon_{\beta_0}$ is very small. In this case, $M_t$
slightly greater than $(A_t+\gamma C_t)/2$ is enough to ensure
exponential decay with rate $\gamma$.

(iv) The case $H\equiv 0$ corresponds to classical MCMC. Here
$A_t=0$ for all $t$, so $\partial\varepsilon_t/\partial t\leq
-2\frac{M_t}{C_t}\varepsilon_t$. This yields the classical
exponential decay with rate $2\gamma$ of the mean square error in
the presence of the global spectral gap $M_t/C_t\geq\gamma$ of the
generator $M_t\cdot {\mathcal L}_t$. For $H\not\equiv 0$,
additional MCMC moves are required to make up for the error growth
due to importance sampling/resampling.
\end{rmk}

Roughly, the corollary says that is the initial error is sufficiently
small, the stabilizing effects of the MCMC dynamics make up for the
error growth due to importance sampling/resampling provided $M_t\geq
A_t/2$.
\medskip

\noindent \emph{Comparison with parallel MCMC.}
Suppose that we want to simulate $\mu_\beta$ for a fixed $\beta>0$.
Parallel MCMC consists in simulating $N$ independent time homogeneous
Markov chains with generator ${\mathcal L}_\beta$. This algorithm is
clearly a special case of the sequential MCMC procedure introduced
above, where $\mu_t=\mu_\beta$ for all $t>0$ and $H=0$. If the chains
are run with initial distribution $\mu_0$, one has
$$
{\varepsilon}_t \ \leq\ e^{-2  t/C_\beta } \,\varepsilon_0\ \leq\
e^{-2 t/C_\beta }\cdot  \big( e^{\beta \osc(H)} - 1 \big)
$$
where we have used that
$$
\varepsilon_0 \ =\  \sum_{x\in S}\left( \frac{\mu_0(x)}{\mu_\beta
(x)} -1\right)^2\,\mu_\beta (x) \ =\ \sum_{x\in S}
\frac{\mu_0(x)}{\mu_\beta(x)}\,\mu_0(x)\, -\,1 \ \leq\ e^{\beta
\osc(H)} - 1\, .
$$
Hence to ensure $\varepsilon_T<\bar\varepsilon$ for a given
$\bar\varepsilon >0$ and $T>0$, a total running time
$$
T \ \geq\ \frac{C_\beta}{2}\cdot\left(\beta\osc(H) +
\log\frac1{\bar\varepsilon} \right)
$$
is sufficient. If (\ref{eq:star}) holds, the number of MCMC steps
required for a simulation is of the same order as $T$.
Alternatively, we can apply the sequential MCMC method with
varying distributions $\mu_t$ ($0\leq t\leq \beta$). Using the
rough estimate $A_t\le C_t\cdot\sup H_t^-$ and (\ref{eq:steps2}),
we see that $\varepsilon_t$ decreases in time if
$$
M_t \ \geq\  \frac 12\,C_t \sup {H}_t^-\,
(1+\varepsilon_0^{1/2})^2 \qquad \forall\ t\in (0,\beta ).
$$
Thus an expected total number of MCMC steps of order
$$
\frac12 (1+\varepsilon_0^{1/2})^2 \int_0^\beta C_t\,\sup H_t^-\, dt
$$
suffices to guarantee stability of the corresponding nonlinear
semigroup.
\bigskip

More drastic improvements due to sequential MCMC appear when good
global spectral gap estimates do not hold, as we shall now
demonstrate.

\subsection{Error control based on local estimates}
\label{subsection:local} Madras and Randall \cite{MR} and Jerrum,
Son, Tetali and Vigoda \cite{JSTV} have shown how to derive
estimates for spectral gaps and logarithmic Sobolev constants of
the generator of a Markov chain from corresponding local estimates
on the sets of a decomposition of the state space combined with
estimates for the projected chain. This has been applied to
tempering algorithms in \cite{MZ}, \cite{BR} and \cite{sch}. We
now develop related decomposition techniques for sequential MCMC.
However, in this case, we will assume {\em only} local estimates
for the generators ${\mathcal L}_t$, and no mixing properties for
the projections -- whence there will be an unavoidable error
growth due to importance sampling/resampling between the
components. The results and examples below indicate that
nevertheless sequential MCMC methods might potentially be at least
equally efficient as tempering algorithms in many
applications.
Since mixing
properties for the projections do not have to be taken into
account, the analysis of the decomposition simplifies
considerably.\medskip

Let $0\leq\beta_0 <\beta_1\leq \infty $. We assume that for every
$t\in (\beta_0,\beta_1)$, there exists a decomposition
$$
S\ =\ \bigcup_{i\in I} S_t^i
$$
into finitely many disjoint sets with $\mu_t(S_t^i)>0$, as well as
non--negative definite quadratic forms ${\mathcal E}_t^i$ ($i\in
I$) on functions on $S$ such that
\begin{equation}
\label{eq:dircomp}
 \sum_i \mu_t(S_t^i)\,\mathcal{E}_t^i(f,f)\
\leq\ K \cdot\mathcal{E}_t(f,f) \qquad\forall \ t\in
(\beta_0,\beta_1),\ f:S\to\r \,
\end{equation}
for some fixed finite constant $K$. For example, one might choose
${\mathcal E}_t^i$ as the Dirichlet form of the Markov chain
corresponding to ${\mathcal L}_t$ restricted to $S_t^i$, i.e.,
\begin{equation}
\label{eq:astast}
{\mathcal E}_t^i(f,f)\ =\ \frac 12\,\sum_{x,y\in
S_t^i}(f(y)-f(x))^2\,{\mathcal L}_t(x,y)\,\mu_t(x\,|\, S_t^i)\, .
\end{equation}
In this case, (\ref{eq:dircomp}) holds with $K=1$.

Let us denote by $\E_t^i$ and $\var_t^i$, respectively, the
expectation and variance w.r.t. the conditional measure
$$
\mu_t^i(A)\ :=\ \mu_t(A|S_t^i),
$$
and by $\pi :S\to I$ the natural projection. In particular,
$$
\E_t[f|\pi]\ =\ \sum_{i\in S}\E_t^i[f]\cdot \chi_{S_t^i},
$$
for any function $f:S\to\mathbb{R}$. We set
$$\tilde{H}_t\ :=\ H -  \E_t[ H|\pi ]\, .$$
Assume that the following {\bf local Poincar\'{e} type
inequalities} hold for all $t\in (\beta_0 ,\beta_1 )$ and $i\in I$
with constants $A_t^i,B_t^i\in (0,\infty )$~:
 \begin{eqnarray}
 \label{eq:Hloc}
\ \E_t^i[-\tilde{H}_t\, f^2] &\leq &A_t^i\cdot
\mathcal{E}_t^i(f,f)
\qquad \forall f:S\to\r\,  :\, \E_t[f|\pi]=0 \ , \\
\label{eq:l1Hloc} \ \big|\E_t^i[\tilde{H}_t\, f]\big|^2 &\leq &
B_t^i\cdot \mathcal{E}_t^i(f,f) \qquad \forall f :S\to\r\,  :\,
\E_t[f|\pi]=0 \ .
\end{eqnarray}
\medskip
\begin{rmk}
\label{rmk:locpoin} (i) Note that to verify (\ref{eq:Hloc}) it is
enough to estimate $\E_t^i[\tilde H_t^-f^2]$, while for
(\ref{eq:l1Hloc}) one has to take into account the positive part
of $\tilde{H}_t$ as well. In particular, (\ref{eq:Hloc}) can not
be used to derive an estimate of type (\ref{eq:l1Hloc}). However,
if (\ref{eq:Hloc}) holds with $-\tilde H_t$ replaced by $|\tilde
H_t|$, then (\ref{eq:l1Hloc}) holds with $B_t^i=\E_t^i[|\tilde
H_t|]\cdot A_t^i$.

(ii) If local Poincar\'e inequalities of the type
 \begin{equation}
 \label{eq:ast}
  \var_t^i(f)\
\leq\  C_t^i\cdot\EE_t^i(f,f)\qquad\forall\ f:S\to\r,\ i\in I,
 \end{equation}
hold, then (\ref{eq:Hloc}) and (\ref{eq:l1Hloc}) hold with
$A_t^i=C_t^i\cdot\max_{S_i}\tilde{H}_t^-$ and
$B_t^i=C_t^i\cdot\var_t^i(H)$.
\end{rmk}

Combining (\ref{eq:dircomp}) and (\ref{eq:Hloc}),
(\ref{eq:l1Hloc}) respectively yields
\begin{equation}\label{eq:LP1}
  \E_t[-\tilde{H}_t\tilde f_t^2]\ =\
      \sum_{i\in I} \mu_t(S_t^i)\, \E_t^i[-\tilde{H}_t \tilde f_t^2]
      \ \leq\  \hat{A}_t \cdot\EE_t(f,f)\qquad\forall\ f:S\to\r\,
      ,
\end{equation}
and
\begin{equation}\label{eq:LP2}
\sum_{i\in I} \mu_t(S_t^i)\, \left|\E_t^i[\tilde{H}_t \tilde
f_t]\right|^2\
      \leq \ \hat{B}_t \cdot\EE_t(f,f)\qquad\ \forall\ f:S\to\r\, ,
 \end{equation}
where $$\hat A_t:=K\cdot\max_i A_t^i\qquad\mbox{and}\qquad \hat
B_t:=K\cdot\max_iB_t^i\, .$$ The following error estimate is our
key result~:

\begin{thm}\label{thm:local}
If $M_t>\hat A_t/2$ for all $t\in (\beta_0,\beta_1)$ then
\begin{equation}
  \label{eq:errloc}
  \frac{d}{dt} \log {\varepsilon}_t \ \leq \ \frac{\hat{B}_t}{M_t-\hat{A}_t/2}\cdot (1+
  {\varepsilon}_t)\,
     + \,\left(1+\sqrt{{\varepsilon}_t}\right)^2\cdot\max_{i\in
     I}h_t^-(i)
\end{equation}
where \begin{equation} \label{eq:ht} h_t(i) \ :=\ \E_t^i[H]\, -\,
\E_t[H]\ =\ -\left.\frac{\partial}{\partial
s}\log\mu_s(S_t^i)\right|_{s=t}\qquad (i\in I) \, .
\end{equation}
\end{thm}\smallskip

The proof is given in Section \ref{sec:proofs} below. To understand the
consequences of (\ref{eq:errloc}), let us first consider the
asymptotics as $M_t$ tends to infinity. In this case,
(\ref{eq:errloc}) reduces to
$$
\frac{d}{dt}\log {\varepsilon}_t \ \leq\ \left(
1+\sqrt{{\varepsilon}_t}\right)^2\cdot \max h_t^-\, .
$$
In order to ensure that for $t>\beta_0$ the error $\varepsilon_t$
remains below a given threshold $\delta>0$, note that as long as
$\varepsilon_t\le\delta$, we have
$$
\frac{d}{dt}\log {\varepsilon}_t \ \leq\ \left(
1+\sqrt\delta\right)^2\cdot\max h_t^-\, .
$$
Thus
\begin{equation}
\label{eq:locgrowth} \min({\varepsilon}_t,\delta )\ \leq\
\varepsilon_{\beta_0}\cdot G_t^{(1+\sqrt\delta )^2} \qquad
\forall\ t\in [\beta_0,\beta_1]\, ,
\end{equation}
where
$$
G_t \ :=\ \exp\left(\int_{\beta_0}^t\max h_r^-\,dr\right)\ =\
\exp\left(\int_{\beta_0}^t\max_i\left.\frac{\partial}{\partial
s}\log\mu_s(S_r^i)\right|_{s=r}\, dr\right).
$$
\begin{rmk}\label{rmk:domi}The term $G_t^{(1+\sqrt\delta )^2}$ in (\ref{eq:locgrowth})
accounts for the maximum error growth due to importance sampling
between the components. If $S_t^i=S^i$ is independent of $t$ for
every $i$, and there is an $i_0\in I$ such that
$\frac{\partial}{\partial s}\log\mu_s(S^i)$ is maximized by
$S^{i_0}$ for all $s\in (\beta_0,\beta_1 )$, then
$$
G_t\ =\ \exp\left( \int_{\beta_0}^t\max_i\frac{d}{d
s}\log\mu_s(S^i)\, ds\right)\ =\ \frac{\mu_t
(S^{i_0})}{\mu_{\beta_0}(S^{i_0})}\qquad\forall\ t\in
[\beta_0,\beta_1],
$$
i.e., $G_t$ is the growth rate of this strongest growing
component. In general, things are more complicated, but a similar
interpretation is at least possible on appropriate subintervals of
$[\beta_0,\beta_1]$.
\end{rmk}

Now we return to the case when $M_t$ is finite. The next corollary
tells us how many MCMC moves are sufficient to obtain an estimate
on the growth of $\varepsilon_t$ that is not much worse than
(\ref{eq:locgrowth}).

\begin{coroll}\label{coroll:locbound}
Let $\beta\in(\beta_0,\beta_1]$ and $\delta >0$, and assume that
\begin{equation}\label{eq:ASS}
M_t \geq \frac{\hat{A}_t}2 + \alpha_t\cdot\hat{B}_t\qquad\forall\
t\in (\beta_0,\beta )
 \end{equation}
 for some function $\alpha:(\beta_0,\beta )\to (0,\infty )$. Then
\begin{equation}
\label{eq:C} \min ({\varepsilon}_\beta ,\delta ) \ \leq \
\varepsilon_{\beta_0}\cdot G_{\beta}^{(1+\sqrt\delta
)^2}\cdot\exp\int_{\beta_0}^\beta\frac{1+\delta }{\alpha_s}\, ds\
. \qquad \forall\ t\in[\beta_0,\beta]\, .
\end{equation}
In particular, if
\begin{equation}\label{eq:BA}
M_t \ \geq\  \frac{\hat{A}_t}2\, +\, (\beta -\beta_0)\cdot\hat
B_t\qquad\forall\ t\in (\beta_0,\beta )
\end{equation}
 then
\begin{equation} \label{eq:CA} \min ({\varepsilon}_\beta ,\delta )
\ \leq \ \varepsilon_{\beta_0}\cdot G_{\beta}^{(1+\sqrt\delta
)^2}\cdot e^{{1+\delta}}\ .
\end{equation}
\end{coroll}\smallskip

\begin{rmk}
The main difference to Corollary \ref{coroll:glob} is that under
local conditions it can not be guaranteed that the error remains
bounded. Instead, $\varepsilon_t$ can grow with a rate dominated
by $G_t^{(1+\sqrt\delta )^2}$. As already pointed out, this is due
to importance sampling between the components.
\end{rmk}

\subsection{Example 1: Exponential model with $k$ valleys in the
energy landscape}
\label{subsection:expmodel} 
This is an extended version of a model considered in \cite{MP},
\cite{MZ} as a test case for some multi-level MCMC methods. We fix
$k\in\n$, and $r_1,r_2,\ldots , r_k\in\n $. Let $S^0:=\{ 0\}$ and
$$
S^i\ :=\ \left\{ (i,j):\, j=1,2,\ldots,r_i\right\}, \qquad 1\leq
i\leq k.
$$
 We consider the graph with vertex set
$$
S \ =\ \bigcup_{i=0}^k S^i
$$
and edges $(0,(i,1))$, $1\leq i \leq k$, and $((i,j),(i,j+1))$,
$1\leq i\leq k$, $1\leq j \leq r_i-1$. Suppose that
$$
H(x)\ =\ -d(x,0),\qquad x\in S,
$$
where $d(x,0)$ stands for the graph distance of $x$ from $0$, i.e.,
$H(0)=0$ and $H((i,j))=-j$. We assume that $\mu_t$ is given by
(\ref{eq:mub}), where $\mu$ is an arbitrary probability distribution
on $S$ such that $\mu(x)>0$ for all $x\in S$ and $\mu$ is log-concave
on each of the valleys $S^i$ of the energy landscape, i.e.,
$$
\frac12 \big( \log\mu((i,j+1))+\log\mu((i,j-1)) \big)
\leq \log\mu((i,j))
$$
for all $1\leq i\leq k$ and $1\leq j\leq r_i$. We consider the
setup for sequential MCMC as described above where $\mathcal{L}_t$
is the generator of the Metropolis dynamics w.r.t. $\mu_t$ based
on the nearest neighbor random walk on $S$. Of course, there are
more efficient ways to carry out Monte Carlo integrations in this
special situation. The point, however, is that sequential MCMC
methods can be applied even though the underlying structure of the
energy landscape is unknown.
Let $R = \max_{1\leq i\leq k} r_i$. An application of Corollary
\ref{coroll:locbound} with $\beta_0=0$ and $S_t^i=S^i$ for all
$t\geq 0$ yields the following result~:

\begin{thm}\label{thm:octopus}
If
$$
M_t \geq R^3+\frac{\beta}2R^4\qquad \forall t\in (0,\beta),
$$
then
\begin{equation}
\label{eq:octo} \min(\varepsilon_\beta,\delta)\ \leq\
e^{1+\delta}\cdot\varepsilon_0
G_\beta^{(1+\sqrt\delta)^2}\cdot\varepsilon_0  \qquad \forall\
\delta\in (0,1).
\end{equation}
Moreover, if the conditional distribution $\mu(\cdot|S^{i_0})$
lies deeper in one of the valleys than in the others in the sense
that
\begin{equation}
  \label{eq:domi}
\mu\big(\{ (i,j):\,j\geq h\}\big| S^{i_0}\big) \ \geq\ \mu\big(\{
(i,j):\,j\geq h\}\big| S^i\big)\ ,
\end{equation}
then
$$
G_\beta\ =\ \frac{\mu_\beta(S^{i_0})}{\mu(S^{i_0})}\ ,
$$
and thus
\begin{equation}
  \label{eq:pus}
\min(\varepsilon_\beta,\delta)\ \leq\ e^{1+\delta}\cdot
\frac{\varepsilon_0}{\mu(S^{i_0})^{(1+\sqrt\delta)^2}}  \qquad
\forall\ 0<\delta <1.
\end{equation}
\end{thm}

\begin{rmk}
  (i) The last estimate indicates that to obtain good bounds it is
  crucial that the mass allocated by the initial distribution on the
  component $S^{i_0}$ with strongest importance growth is not too
  small (although it can be rather small if the initial distribution
  $\nu_0$ is a good approximation of $\mu_0$).

  (ii) Let $K_\beta=\int_0^\beta M_t\,dt$. Note that $K_\beta$ is a
  measure for the total number of MCMC steps that a corresponding
  sequential MCMC algorithm will perform on average. The theorem
  implies that choosing $M_t$ constant on $[0,\beta]$ with $K_\beta$
  of order $O(\beta^2)$ is sufficient to guarantee that the nonlinear
  flow of measures has good stability properties on $[0,\beta]$, and
  can thus be used to efficiently approximate $\mu_\beta$. In contrast
  to this situation, the flow of measures corresponding to the
  standard
  simulated annealing algorithm has good stability properties only
  if $K_\beta$ grows exponentially in $\beta$.
\end{rmk}

\subsection{Example 2: The mean field Ising model}
\label{subsection:ising}
As a very simple example for a model with
a phase transition, we now consider the mean field Ising
(Curie--Weiss) model, i.e. $\mu_\beta$ is of type (\ref{eq:mub})
where $\mu_0=\mu $ is the uniform distribution on the hypercube
$$S\ =\ \{ -1,+1\}^N\, ,$$
and
\begin{equation}
H(\sigma )\ =\ -\frac
1{2N}\,\sum_{i,j=1}^N \sigma_i\sigma_j
\end{equation}
for some $N\in\n $. Let ${\mathcal L}_\beta $ be the generator of
the (time--continuous) Metropolis chain w.r.t.\ $\mu_\beta $ based
on the nearest neighbor random walk on $S$ as proposal matrix. It
is well known that this chain is rapidly mixing (i.e.\ the
spectral gap decays polynomially in $N$) for $\beta < 1$, but
torpid mixing holds (i.e.\ the spectral gap decays exponentially
fast in $N$) for $\beta >1$. Thus in the multi-phase regime $\beta
>1$, the classical Metropolis algorithm converges to equilibrium
extremely slowly for large $N$.

Now assume for simplicity that $N$ is odd, and decompose $S$ into
the two components
\begin{eqnarray*}
S^+ & :=& \left\{ \sigma\in S\, \left|\ \sum_{i=1}^N\sigma_i >
0\right.\right\} \qquad\mbox{and}\\
S^- & :=& \left\{ \sigma\in S\, \left|\ \sum_{i=1}^N\sigma_i <
0\right.\right\} .
\end{eqnarray*}
Improving on previous results (e.g. of Madras and Zheng
\cite{MZ}), Schweizer \cite{sch} showed recently that the spectral
gaps of the restricted Metropolis chains on both $S^+$ and $S^-$
are bounded from below by $\frac{1}{9}N^{-2}$ for {\em every
}$t\ge 0$.
Applying the results above to the error growth for the non-linear
semigroup corresponding to sequential MCMC in this situation, we
obtain~:

\begin{thm}\label{cor:ising}
For every $\beta >0$ and $N\in \n $,
$$\sup_{0\le t\le\beta }\, \varepsilon_t\ \le\
 e^2\cdot\varepsilon_0$$ holds whenever $\varepsilon_0\le 1$ and
\begin{equation} \label{eq:ss}
M_t\ \ge \ \frac 94\, N^3\, +\,\frac 98\,\beta\, N^4\qquad\forall\
t\in (0,\beta ) .
\end{equation}
\end{thm}

\begin{rmk}
(i) The result is based on a rough estimate of $\hat A_t$ and
$\hat B_t$ in terms of the local spectral gap. We expect that a
more precise estimate of these constants would yield a smaller
power of $N$ in (\ref{eq:ss}). Furthermore, for $\beta\leq 1$, the
result can be improved by applying global instead of local
spectral gap estimates. However, our main interest is the phase
transition regime.

(ii) Related results for the mean field Ising model have been
obtained for mixing times of Markov chains for umbrella sampling
in \cite{MP}, and for simulated and parallel tempering in
\cite{MZ}, \cite{BR}, \cite{sch}. Schweizer \cite{sch} obtains an
upper bound on the order in $N$ and $\beta$ of the $L^2$ mixing
time (inverse spectral gap) for simulated tempering that is close
to the one in (\ref{eq:ss}). In contrast, the best known order for
parallel tempering is much worse. In general, it seems that the
analysis of sequential MCMC is partially simpler than the one for
parallel tempering, where one has to take into account that a
particle can only move in temperature if another particle moves in
the opposite direction.  In fact, for this reason we would expect
that sequential MCMC methods can have substantial advantages
compared to parallel tempering.

(iii) The theorem can be extended to a mean field Ising model with
magnetic field. In this case, however, one has to take into
account an additional (but well controlled) error growth due to
importance sampling/resampling between the components. Moreover,
the decomposition into the two components will now depend on $t$.
Without magnetic field this is not the case because of the
built-in symmetry.
\end{rmk}

\subsection{Extensions}
\label{subsection:extensions} As remarked above, our results
immediately extend to the case where
$$
 \mu_t(x)\ =\ \frac 1{Z_t}\,e^{-\int_0^tU_s(x)\,ds}\,\mu (x)\,
 \qquad (0\leq t\leq\beta )
 $$
 for a continuous function $(s,x)\mapsto U_s(x)$ on $[0,\beta ]$
 with $U_s(x)\geq 0$ for all $s\in [0,\beta ]$ and $x\in S$. In
 this case, the evolution equations (\ref{eq:nut}) and
 (\ref{eq:X}) in Theorem~1 take the form
\begin{eqnarray*}
\frac{\partial}{\partial t} {\nu}_t&  =& M_t\nu_t{\mathcal L}_t\,
-\, U_t\,\nu_t\, +\, \nu_t(U_t)\, \nu_t, \qquad \mbox{and}\\
\frac{\partial}{\partial t} {g}_t&  =& M_t{\mathcal L}_tg_t\, +\,
\E_t[U_t (g_t-1)]\, g_t.
\end{eqnarray*}
The evolution equation (\ref{eq:err}) for the mean square error
and all the stability estimates in Sections
\ref{subsection:global} and \ref{subsection:local} still hold if
$H$ is replaced by $U_t$, and, correspondingly,
$$H_t\ =\ U_t\, -\, \E_t[U_t]\, .$$
All proofs are completely analogous.\medskip

The extension of the results to more general state spaces requires
some (standard) technical assumptions which make the proofs
slightly less transparent. We postpone this extension to a future
publication where we will also consider corresponding
applications.

\section{Proofs}
\label{sec:proofs}
\subsection{Proof of Theorem \ref{thm:err}}
To simplify the notation, we assume $M_t=1$ for all $t\geq 0$. The
general case is similar with ${\mathcal L}_t$ replaced by
$M_t\cdot {\mathcal L}_t$. Let us also set $p_t:= p_{0,t}$ and
$\Phi_t:=\Phi_{0,t}$. Then one has
$$
\nu_t\ =\ \Phi_t\nu\ =\ \frac{\nu p_t}{(\nu p_t)(1)}\, .
$$
The forward equation (\ref{eq:fp}) yields
$$
\frac{\partial}{\partial t}\nu p_t\  =\ \nu p_t {\mathcal L}_t
\,-\, H\,\nu p_t\, .
$$
Since $(\nu p_t)(1)>0$ and ${\mathcal L}_t1=0$ for all $t$, we
obtain
\begin{eqnarray}
\nonumber\frac{\partial}{\partial t} \nu_t &=&
\frac{\partial}{\partial t}\, \frac{\nu p_t}{{(\nu p_t)}(1)} \\
& =& \frac{\nu p_t{\mathcal L}_t\, -\, H\nu p_t}{(\nu p_t)(1)}\,
-\,\frac{(\nu p_t {\mathcal L}_t-H\nu p_t)(1)\,\nu p_t}{(\nu
p_t)(1)^2}
 \\
 \label{eq:q}
&=& \nu_t{\mathcal L}_t\, -\, H\nu_t\, +\,\nu_t(H)\, \nu_t\ \, .
\end{eqnarray}

Next, we derive a corresponding evolution equation for the
densities
$$
g_t(y)\  :=\ \frac{\nu_t(y)}{\mu_t(y)}\qquad (y\in S).
$$
Since $\mu_t$ has full support and is differentiable in $t$, we
obtain
\begin{equation}
\label{eq:h1} \frac{\partial}{\partial t} g_t\ =\ \frac
1{\mu_t}\,\frac{\partial}{\partial t} \nu_t\,-\,
\frac{\nu_t}{\mu_t}\, \frac{\partial}{\partial t} \log\mu_t\ .
\end{equation}
Note that by the detailed balance condition (\ref{eq:db}), the
relative density of $\nu_t{\mathcal L}_t$ w.r.t.\ $\mu_t$ is
$$\frac{(\nu_t{\mathcal L}_t)(y)}{\mu_t(y)}\ =\
\sum_x\nu_t(x)\,\frac{{\mathcal L}_t(x,y)}{\mu_t(y)}\ =\
\sum_x\nu_t(x)\,\frac{{\mathcal L}_t(y,x)}{\mu_t(x)}\ =\
({\mathcal L}_tg_t)(y)\, .$$
 Hence (\ref{eq:q})
yields
\begin{eqnarray}
\frac 1{\mu_t}\, \frac{\partial}{\partial t}\,\nu_t  &=&
    {\mathcal L}_tg_t-Hg_t+\nu_t(H)\, g_t\nonumber\\
    &=& ({\mathcal L}_t-H)\, g_t\, +\, \E_t[Hg_t]\, g_t\, . \label{eq:pezzo1}
\end{eqnarray}
Recalling that $\mu_t = \frac1{Z_t}e^{-tH}\mu$ with $Z_t=\sum_S
e^{-tH(y)}\,\mu (y)$, one has
\begin{equation}
  \label{eq:mu}
  \frac{\partial}{\partial t} \log\mu_t\ =\ -\mu_t \left( H - \E_t[H] \right)\ =\ -\mu_tH_t\ ,
\end{equation}
hence
\begin{equation}
\label{eq:pezzo2}
  -\frac{\nu_t}{\mu_t}\, \frac{\partial}{\partial t} \log\mu_t\ =\
  \left( H-\E_t[H]\right)\, g_t\ .
\end{equation}
Inserting (\ref{eq:pezzo1}) and (\ref{eq:pezzo2}) into
(\ref{eq:h1}) we obtain
\begin{eqnarray}
\frac{\partial}{\partial t} g_t &=&
    {\mathcal L}_tg_t\, +\,  \big( \E_t[Hg_t] - \E_t[H] \big)\, g_t \nonumber\\
&=& {\mathcal L}_tg_t \, +\,  \E_t\big[ H(g_t-1) \big]\, g_t\, .
\label{eq:h2}
\end{eqnarray}
We are now ready to derive the equation for the quadratic error
$$
  {\varepsilon}_t \ =\ \E_t\big[ (g_t-1)^2 \big].
$$
Differentiating this expression with respect to $t$ one gets by
(\ref{eq:h2}) and (\ref{eq:mu}),
\begin{eqnarray}
  \nonumber
  \lefteqn{\frac{d}{dt} {\varepsilon}_t \ =\ 2\, \E_t\left[(\frac{\partial}{\partial t} g_t)(g_t-1)\right]
  \,+\,           \E_t\left[ (g_t-1)^2\,\frac{\partial}{\partial t}\log\mu_t
  \right] }\\
  \nonumber
  &=&
2\,\E_t\left[ ({\mathcal L}_tg_t)(g_t-1) \right]\, +\,
   2\,\E_t\left[ g_t(g_t-1) \right] \, \E_t\left[ H(g_t-1) \right]\, -\,
   \E_t\left[ H_t (g_t-1)^2 \right]\\
  &=& 2\,\E_t\left[ {\mathcal L}_t(g_t-1)\, (g_t-1) \right]
\, +\,
   2\,\E_t\left[ (g_t-1)^2 \right] \, \E_t\left[ H(g_t-1) \right]
       \, -\, \E_t\left[ H_t(g_t-1)^2 \right]  \nonumber\\
     &=& -2\, \mathcal{E}_t(g_t-1) \, +\,
   2\, \E_t\left[ H(g_t-1) \right]\cdot\varepsilon_t\, -\, \E_t\left[ H_t(g_t-1)^2
   \right]\ .
 \nonumber
\end{eqnarray}
In the derivation we have used that $$\E_t[g_t-1]\ =\
\nu_t(1)-\mu_t(1)\ =\ 0\, ,$$ and ${\mathcal L}_t1\equiv 0$. The
equation implies (\ref{eq:err}) in the case $M_t\equiv 1$. The
general case follows similarly. \hfill$\square$

\subsection{Proof of Theorem \ref{thm:logerr}}

We have to estimate the terms on the right hand side of
(\ref{eq:err}). By the assumed $H$--Poincar\'{e} inequality
(\ref{eq:H}), we obtain
$$
-\frac12\, \E_t\left[H_t(g_t-1)^2\right]  \ \leq\ \frac12
\,\E_t\left[H_t^-(g_t-1)^2\right]
      \ \leq \ \frac12 A_t\cdot \EE_t(g_t-1)\ .
$$
Moreover,
$$
\E_t\left[H_t(g_t-1)\right] \ \leq\ \left(\E_t[H_t^2]\right)^{1/2}
\left(\E_t[(g_t-1)^2]\right)^{1/2}\ =\ \sigma_t(H)\,
{\varepsilon}_t^{1/2}.
$$
Substituting into (\ref{eq:err}) yields
\begin{eqnarray*}
\frac{d}{dt} {\varepsilon}_t &\leq& -2\, (M_t-A_t/2)\,\EE_t(g_t-1) \, +\, 2\,\sigma_t(H)\, {\varepsilon}_t^{3/2} \\
&\leq& -\frac{2M_t-A_t}{C_t}\,{\varepsilon}_t \, +\,
2\,\sigma_t(H)\, {\varepsilon}_t^{3/2},
\end{eqnarray*}
by the global Poincar\'e inequality (\ref{eq:Poin}), provided
$M_t\geq A_t/2$. This proves (\ref{eq:log1}).
\medskip

On the other hand,
\begin{eqnarray}
\nonumber
\lefteqn{\E_t\left[ H_t\big(-(g_t-1)^2/2+(g_t-1){\varepsilon}_t\big)\right]}  \\
\label{eq:A1} &=& \frac12 \,\E_t\left[ H_t^-(g_t-1)^2 \right] +
\E_t\left[ H_t^-(1-g_t) \right]\,  {\varepsilon}_t  \\
\nonumber &&\qquad + \E_t\left[
H_t^+\big(-(g_t-1)^2/2+(g_t-1){\varepsilon}_t\big)\right] \, .
\end{eqnarray}
Estimating the three summands on the right hand side separately
yields
$$
\E_t\big[ H_t^-(g_t-1)^2 \big] \ \leq\ A_t\cdot\EE_t(g_t-1)
$$
by the $H$-Poincar\'e inequality (\ref{eq:H}),
\begin{eqnarray*}
\E_t\big[ H_t^-(1-g_t) \big]
&\leq& \E_t[H_t^-]^{1/2} \, \E_t[H_t^-(g_t-1)^2]^{1/2} \\
&\leq& \E_t[H_t^-]^{1/2} A_t^{1/2} \EE_t(g_t-1)^{1/2}
\end{eqnarray*}
by the Cauchy-Schwarz inequality and (\ref{eq:H}), and
$$
\E_t\big[ H_t^+\big(-(g_t-1)^2/2+(g_t-1){\varepsilon}_t\big)\big]
\ \leq\ \frac 12\,\E_t[H_t^+] \, {\varepsilon}_t^2 \ =\ \frac 12
\,\E_t[H_t^-] \, {\varepsilon}_t^2\ .
$$
The last estimate follows since
$$\frac12 \, {\varepsilon}_t^2\ \geq\
(g_t-1){\varepsilon}_t - \frac12 (g_t-1)^2$$ and
$$\E_t[H_t^+]-\E_t[H_t^-]\ =\ \E_t[H_t]\ =\ 0\ .$$ By combining the
estimates, (\ref{eq:A1}) and (\ref{eq:err}), we obtain
$$
\frac{d}{dt} {\varepsilon}_t \ \leq\ -(2M_t-A_t)\,\EE_t(g_t-1) \,
+\, 2 A_t^{1/2} \E_t[H_t^-]^{1/2} \EE_t(g_t-1)^{1/2}
{\varepsilon}_t \,  +\, \E_t[H_t^-] \, {\varepsilon}_t^2\, .
$$
This combined with the global Poincar\'e inequality
(\ref{eq:Poin}) yields
$$
\frac{d}{dt} {\varepsilon}_t\ \leq\ -\frac{2M_t-A_t}{C_t}\,
{\varepsilon}_t \, + \, 2\, \frac{A_t^{1/2}}{C_t^{1/2}}\,
\E_t[H_t^-]^{1/2} {\varepsilon}_t^{3/2}\, +\, \E_t[H_t^-]\,
{\varepsilon}_t^2\, ,
$$
and hence (\ref{eq:log2}). \hfill$\square$

\subsection{Proof of Corollary \ref{coroll:glob}}

If (\ref{eq:steps1}) or (\ref{eq:steps2}) holds for $t\in
(\beta_0,\beta_1)$, then by Theorem \ref{thm:logerr} and
continuity, $t\mapsto\varepsilon_t$ is strictly decreasing near
$\beta_0$ and near any $s\in (\beta_0,\beta_1 )$ such that
$\varepsilon_s\leq\varepsilon_{\beta_0}$. Hence it is strictly
decreasing on the whole interval $[\beta_0
,\beta_1]$.\hfill$\square$

\subsection{Proof of Theorem \ref{thm:local}}

Similarly to Theorem \ref{thm:logerr}, we have to control the
right hand side of (\ref{eq:err}), but now by using only local
Poincar\'{e} type inequalities. Let
$$f_t\ :=\ g_t-1\qquad\mbox{and}\qquad
\tilde{f}_t\ :=\ f_t-\E_t[f_t|\pi ]\ .$$ Then
\begin{eqnarray}
\nonumber \lefteqn{\E_t\left[
H_t\cdot\big(-(g_t-1)^2/2+(g_t-1){\varepsilon}_t\big) \right] }\\
\nonumber
&=& \E_t\left[ \tilde{H}_t\big(-(g_t-1)^2/2+(g_t-1){\varepsilon}_t\big) \right] \\
\nonumber &&    +\, \sum_{i\in
I}\mu_t(S_t^i)\left(\E_t^i[H]-\E_t[H]\right)\cdot
        \E_t^i\left[-(g_t-1)^2/2+(g_t-1){\varepsilon}_t\right]\\
\nonumber &=&  -\frac12\, \E_t[\tilde{H}_t\tilde{f}_t^2]\,
-\,\E_t\big[\tilde{H}_t\tilde{f}_t\,\E_t[ f_t|\pi ]\big]
\, +\,\E_t[\tilde{H}_t\tilde{f}_t\,{\varepsilon}_t]\\
\label{eq:AAA}
 && \,  + \sum_{i\in I}\mu_t(S_t^i)\cdot
 \left(\E_t^i[H]-\E_t[H]\right)\cdot
        \E_t^i\left[-f_t^2/2+f_t{\varepsilon}_t\right] \\
\nonumber
 &=&  -\frac12\, \E_t[\tilde{H}_t\tilde{f}_t^2]\,
 +\, \sum_{i\in I} \mu_t(S_t^i)\cdot \E_t^i[\tilde{H}_t\tilde{f}_t]\cdot\left({\varepsilon}_t-\E_t^i[f_t]\right) \\
\nonumber && \, + \sum_{i\in I} \mu_t(S_t^i)\,
h_t(i)\cdot\E_t^i[-f_t^2/2+f_t{\varepsilon}_t]\, .
\end{eqnarray}
Here we have used the definitions of $H_t$, $\tilde H_t$ and
$h_t$, and the fact that $\E_t[\tilde H_t|\pi ]=0$. We now
estimate the three summands on the right hand side separately. By
the local $H$-Poincar\'{e} inequality (\ref{eq:LP1}),
$$
-\frac12 \, \E_t [\tilde{H}_t\tilde{f}_t^2]\ \leq \frac12\,
\hat{A}_t\cdot \EE_t(f_t).
$$
By (\ref{eq:LP2}), and since
$$\sum_i\mu_t(S_t^i)\,\E_t^i[f_t]\ =\ \E_t[f_t]\ =\ 0\ ,$$
we have
\begin{eqnarray*}
\lefteqn{\sum_{i\in I} \mu_t(S_t^i)\cdot
\E_t^i[\tilde{H}_t\tilde{f}_t]\cdot
\left({\varepsilon}_t-\E_t^i[f_t]\right)}\\
&\leq& \left(\sum_{i\in I} \mu_t(S_t^i)\,\E_t^i[\tilde H_t\tilde
f_t]^2\right)^{1/2}\left(
             \sum_{i\in I} \mu_t(S_t^i)\, ({\varepsilon}_t-\E_t^i[f_t])^2\right)^{1/2} \\
&\leq& \hat{B}_t^{1/2} \EE_t(f_t)^{1/2}\cdot \left(
\varepsilon_t^2+\sum\mu_t(S_t^i)\E_t^i[f_t^2]\right)^{1/2}\\
&=& \left( \hat B_t\,\EE_t(f_t)\cdot {\varepsilon}_t\cdot
(1+\varepsilon_t)\right)^{1/2}\, .
\end{eqnarray*}
Moreover, since
$$-f_t^2/2+f_t\varepsilon_t\ \leq\ \varepsilon_t^2/2\ ,$$ we obtain
\begin{eqnarray*}
\lefteqn{\sum_{i\in I} \mu_t(S_t^i)\, h_t(i)\,
\E_t^i[-f_t^2/2+f_t{\varepsilon}_t]}\\
&\leq & \sum_{i\in I} \mu_t(S_t^i)\, h_t^+(i)\cdot\frac
12\,\varepsilon_t^2\, +\,\sum_{i\in I} \mu_t(S_t^i)\, h_t^-(i)\,
\E_t^i[f_t^2/2-f_t{\varepsilon}_t]\\
&\leq & \left( \frac 12\,\varepsilon_t^2\, +\, \frac
12\,\varepsilon_t\, +\,\varepsilon_t^{3/2}\right)\cdot\max h_t^-\
=\ \varepsilon_t\cdot \left(
1+\sqrt{\varepsilon_t}\right)^2\cdot\max h_t^-\, .
\end{eqnarray*}
Here we have used that
$$\sum\mu_t(S_t^i)h_t^+(i)\ =\ \sum\mu_t(S_t^i)h_t^-(i)\ \leq\
\max h_t^-\,
,\qquad\mbox{and}$$
$$\sum\mu_t(S_t^i)\,\E_t^i[-f_t]\ \leq\ \left(\sum\mu_t(S_t^i)\, \E_t^i[f_t^2]\right)^{1/2}\ =\
{\varepsilon}_t^{1/2}\, .
$$\smallskip

Combining the estimates yields by (\ref{eq:err}) and
(\ref{eq:AAA})~:
\begin{eqnarray*}
\lefteqn{\frac 12 \frac{d}{dt} {\varepsilon}_t \ \leq\ -M_t\cdot
\EE_t(f_t)\,
+\, \E_t\big[ H_t(-f_t^2/2 + f_t{\varepsilon}_t) \big] } \\
&\leq& -\left(M_t-\frac{\hat{A}_t}{2}\right)\cdot\EE_t(f_t)\, +\,
\left(\hat{B}_t\EE_t(f_t){\varepsilon}_t(1+\varepsilon_t)\right)^{1/2}\,
+\, \frac12 \, {\varepsilon}_t(1+\sqrt{{\varepsilon}_t})^2\max h_t^- \\
&\leq& \frac{\hat{B}_t}{2M_t-\hat{A}_t}\, \varepsilon_t\,
(1+\varepsilon_t)\, +\, \frac12 \, {\varepsilon}_t\, \left(
1+\sqrt{{\varepsilon}_t}\right)^2\max h_t^-\ ,
\end{eqnarray*}
provided $M_t>\hat A_t/2$. This proves
(\ref{eq:errloc}).\smallskip

Moreover, for any subset $A\subseteq S$,
 \begin{eqnarray*}
 \frac{d}{dt}\log\mu_t(A) &=& \frac{d}{dt}\log\sum_{x\in
 A}e^{-tH(x)}\mu (x)\, -\, \frac{d}{dt}\log Z_t\\
 &=& -\E_t[H|A]\, +\, \E_t[H]\, ,
 \end{eqnarray*}
 which proves (\ref{eq:ht}).
 \hfill$\square$

\subsection{Proof of Corollary \ref{coroll:locbound}}

Assume that (\ref{eq:ASS}) holds, and let
$$u_t\ :=\ {\varepsilon}_t/G^{(1+\sqrt\delta )^2}_t\ .$$
Then by the
definition of $G_t$, Theorem \ref{thm:local}, and (\ref{eq:ASS}),
\begin{eqnarray*}
\frac{d}{dt}\log u_t &=& \frac{d}{dt} \log {\varepsilon}_t\, -\,
(1+\sqrt\delta )^2\max h_t^-\\  &\leq & \frac{\hat B_t}{M_t-\hat
A_t/2}\cdot (1+\delta )\  \leq\ \frac{1+\delta}{\alpha_t}
\end{eqnarray*}
for all $t\in (\beta_0,\beta)$ such that
${\varepsilon}_t\leq\delta$. Hence
$$
\varepsilon_t\ =\ u_t\cdot G_t^{(1+\sqrt\delta )^2}\ \leq\
\varepsilon_{\beta_0}\cdot\exp\int_{\beta_0}^t\frac{1+\delta}{\alpha_s}\,
ds\cdot G_t^{(1+\sqrt\delta )^2}$$ holds for $t\in [\beta_0,\beta
]$ provided the right hand side is smaller than $\delta$. This
proves (\ref{eq:C}). The second assertion is a straightforward
consequence.\hfill$\square$

\subsection{Proof of Theorem \ref{thm:octopus}}

The log-concavity of $\mu$ easily implies that $\mu_t$ as well is
log-concave on $S^i$ for all $t\geq 0$ and $1\leq i\leq k$. In
particular, the restriction of $\mu_t$ to $S^i$ has a unique local
maximum for every $i$. By the method of canonical paths it is then
not difficult to prove that the spectral gap of the Metropolis
dynamics w.r.t. $\mu_t(\cdot|S^i)$ based on the standard random
walk is bounded from below by $1/2r_i^2$ for all $t\geq 0$ and
$1\leq i\leq k$, cf.\ e.g.\ Proposition 6.3 in \cite{DSC-metro}.
Now we are in the setting of Remark \ref{rmk:locpoin} (ii),
according to which (\ref{eq:Hloc}) and (\ref{eq:l1Hloc}) hold with
$\EE_t^i$ as in (\ref{eq:astast}),
$$
A_t^i = 2r_i^3\, ,  \qquad\mbox{and}\qquad B_t^i = \frac12 r_i^4\
.
$$
Estimate (\ref{eq:octo}) now follows by a straightforward application
of Corollary \ref{coroll:locbound}.

To prove the second part of the assertion, we show that
(\ref{eq:domi}) places us in the setting of Remark \ref{rmk:domi}.
In fact, for $t>0$,
$$\frac{d}{dt}\log\mu_t(S^i)\ =\ \E_t[H]-\E_t^i[H]\qquad\mbox{for
all }i, \ \mbox{ and}$$
$$-\E_t^i[H]\ =\ -\frac{\mu \left( \left. He^{-tH}\right| S_i\right) }{\mu \left( \left. e^{-tH}\right| S_i\right)
}\ =\ \frac{\sum_jje^{tj}\mu ((i,j))}{\sum_je^{tj}\mu ((i,j))}\
.$$
 If (\ref{eq:domi}) holds, then for any $t>0$, the right hand
side is maximized when $i=i_0$. Hence by Remark \ref{rmk:domi},
$$G_t\ =\ \frac{\mu_t(S_{i_0})}{\mu (S_{i_0})} \qquad\mbox{for all
} t\geq 0.$$ \hfill$\square$

\subsection{Proof of Theorem \ref{cor:ising}}

Since $-N/2\leq H(\sigma )\leq 0$ for all $\sigma $, we have ${\rm
osc}\,(H)\leq N/2$ and
$${\rm Var}_t(H|S^+)\ =\ {\rm Var}_t(H|S^-)\ \leq\ \left(\frac 12\, {\rm
osc}\, (H)\right)^2\ \leq\ N^2/8$$ for every $t\geq 0$. By
Schweizer's result \cite{sch}, a local Poincar\'{e} inequality of
type (\ref{eq:ast}) holds on $S^+$ and $S^-$ with
$C_t^+=C_t^-=9N^2$. Hence by Remark \ref{rmk:locpoin}~(ii),
(\ref{eq:Hloc}) and (\ref{eq:l1Hloc}) hold with
$$A_t^\pm \ =\ \frac 92\, N^3\qquad\mbox{and}\qquad B_t^\pm\ =\
\frac 98\, N^4\, .$$ The assertion now follows from Corollary
\ref{coroll:locbound}, since
$$\E_t^+[H]\ =\ \E_t^-[H]\ =\ \E_t[H]\ .$$\hfill$\square$

\bibliographystyle{amsplain}

\begin{thebibliography}{10}

\bibitem{BR}
N.~Bhatnagar and D.~Randall, \emph{Torpid mixing of simulated tempering on the
  {P}otts model}, {SODA} '04: {P}roceedings of the fifteenth annual {ACM-SIAM}
  symposium on Discrete algorithms (Philadelphia, PA, USA), Society for
  Industrial and Applied Mathematics, 2004, pp.~478--487.

\bibitem{BG}
R.~M. Blumenthal and R.~K. Getoor, \emph{Markov processes and potential
  theory}, Pure and Applied Mathematics, Vol. 29, Academic Press, New York,
  1968. \MR{MR0264757 (41 \#9348)}

\bibitem{cappe}
O.~Capp{\'e}, E.~Moulines, and T.~Ryd{\'e}n, \emph{Inference in hidden {M}arkov
  models}, Springer Series in Statistics, Springer, New York, 2005.
  \MR{MR2159833 (2006e:60002)}

\bibitem{DM}
P.~Del~Moral, \emph{Feynman-{K}ac formulae}, Springer-Verlag, New York, 2004.
  \MR{MR2044973 (2005f:60003)}

\bibitem{DMD}
P.~Del~Moral and A.~Doucet, \emph{On a class of genealogical and interacting
  {M}etropolis models}, S\'eminaire de Probabilit\'es XXXVII, Lecture Notes in
  Math., vol. 1832, Springer, Berlin, 2003, pp.~415--446. \MR{MR2053058
  (2005g:65013)}

\bibitem{DMDJ}
P.~Del~Moral, A.~Doucet, and A.~Jasra, \emph{Sequential {M}onte {C}arlo
  samplers}, J. R. Statist. Soc. B \textbf{68} (2006), no.~3, 411--436.
  \MR{MR1819122 (2002k:60013)}

\bibitem{DMG}
P.~Del~Moral and A.~Guionnet, \emph{On the stability of interacting processes
  with applications to filtering and genetic algorithms}, Ann. Inst. H.
  Poincar\'e Probab. Statist. \textbf{37} (2001), no.~2, 155--194.
  \MR{MR1819122 (2002k:60013)}

\bibitem{MM}
P.~Del~Moral and L.~Miclo, \emph{Branching and interacting particle systems
  approximations of {F}eynman-{K}ac formulae with applications to non-linear
  filtering}, S\'eminaire de Probabilit\'es, XXXIV, Lecture Notes in Math.,
  vol. 1729, Springer, Berlin, 2000, pp.~1--145. \MR{MR1768060 (2001g:60091)}

\bibitem{DSC-comp}
P.~Diaconis and L.~Saloff-Coste, \emph{Comparison theorems for reversible
  {M}arkov chains}, Ann. Appl. Probab. \textbf{3} (1993), no.~3, 696--730.
  \MR{MR1233621 (94i:60074)}

\bibitem{DSC-logsob}
\bysame, \emph{Logarithmic {S}obolev inequalities for finite {M}arkov chains},
  Ann. Appl. Probab. \textbf{6} (1996), no.~3, 695--750. \MR{MR1410112
  (97k:60176)}

\bibitem{DSC-nash}
\bysame, \emph{Nash inequalities for finite {M}arkov chains}, J. Theoret.
  Probab. \textbf{9} (1996), no.~2, 459--510. \MR{MR1385408 (97d:60114)}

\bibitem{DSC-metro}
\bysame, \emph{What do we know about the {M}etropolis algorithm?}, J. Comput.
  System Sci. \textbf{57} (1998), no.~1, 20--36, 27th Annual ACM Symposium on
  the Theory of Computing (STOC'95) (Las Vegas, NV). \MR{MR1649805
  (2000b:68094)}

\bibitem{EM2}
A.~Eberle and C.~Marinelli, \emph{Convergence of sequential {M}arkov chain
  {M}onte {C}arlo methods: {II}. {A}symptotic analysis of interacting particle
  systems}, In preparation.

\bibitem{JSTV}
M.~Jerrum, J.-B. Son, P.~Tetali, and E.~Vigoda, \emph{Elementary bounds on
  {P}oincar\'e and log-{S}obolev constants for decomposable {M}arkov chains},
  Ann. Appl. Probab. \textbf{14} (2004), no.~4, 1741--1765. \MR{MR2099650
  (2005i:60139)}

\bibitem{liu}
J.~S. Liu, \emph{Monte {C}arlo strategies in scientific computing},
  Springer-Verlag, New York, 2001. \MR{MR1842342 (2002i:65006)}

\bibitem{MP}
N.~Madras and M.~Piccioni, \emph{Importance sampling for families of
  distributions}, Ann. Appl. Probab. \textbf{9} (1999), no.~4, 1202--1225.
  \MR{MR1728560 (2001e:60139)}

\bibitem{MR}
N.~Madras and D.~Randall, \emph{Markov chain decomposition for convergence rate
  analysis}, Ann. Appl. Probab. \textbf{12} (2002), no.~2, 581--606.
  \MR{MR1910641 (2003d:60135)}

\bibitem{MZ}
N.~Madras and Z.~Zheng, \emph{On the swapping algorithm}, Random Structures
  Algorithms \textbf{22} (2003), no.~1, 66--97. \MR{MR1943860 (2004c:82117)}

\bibitem{casella}
C.~P. Robert and G.~Casella, \emph{Monte {C}arlo statistical methods}, second
  ed., Springer-Verlag, New York, 2004. \MR{MR2080278 (2005d:62006)}

\bibitem{SC}
L.~Saloff-Coste, \emph{Lectures on finite {M}arkov chains}, Lectures on
  probability theory and statistics (Saint-Flour, 1996), Lecture Notes in
  Math., vol. 1665, Springer, Berlin, 1997, pp.~301--413. \MR{MR1490046
  (99b:60119)}

\bibitem{sch}
N.~Schweizer, \emph{Diploma thesis, {U}niversit\"at {B}onn}, 2006.

\bibitem{St-var}
W.~Stannat, \emph{On the convergence of genetic algorithms---a variational
  approach}, Probab. Theory Related Fields \textbf{129} (2004), no.~1,
  113--132. \MR{MR2052865 (2005d:35040)}

\end{thebibliography}

\def\cprime{$'$}
\providecommand{\bysame}{\leavevmode\hbox to3em{\hrulefill}\thinspace}
\providecommand{\MR}{\relax\ifhmode\unskip\space\fi MR }
\providecommand{\MRhref}[2]{%
  \href{http://www.ams.org/mathscinet-getitem?mr=#1}{#2}
}
\providecommand{\href}[2]{#2}

\end{document}